\documentclass[a4paper,11pt,twoside,reqno]{amsart}

\usepackage[utf8]{inputenc}
\usepackage[plainpages=false,pdfpagelabels=true]{hyperref}
\usepackage{amsmath}
\usepackage{amssymb,amsthm}
\usepackage[margin=1in]{geometry}
\usepackage{slashed}

\newtheorem{Satz}{Theorem}[section]
\newtheorem{Prop}[Satz]{Proposition}
\newtheorem{Lem}[Satz]{Lemma}

\newtheorem{Cor}[Satz]{Corollary}
\newtheorem{Con}[Satz]{Conjecture}
\theoremstyle{definition}
\newtheorem{Dfn}[Satz]{Definition}
\newtheorem{Bem}[Satz]{Remark}

\newcommand{\s}{{\mathbb{S}}}
\allowdisplaybreaks[1]

\renewcommand{\epsilon}{\varepsilon}

\newcommand{\R}{\ensuremath{\mathbb{R}}}

\numberwithin{equation}{section}

\title{Classification results for polyharmonic helices in space forms}
\author{Volker Branding}
\date{\today}
\address{University of Vienna, Faculty of Mathematics\\
Oskar-Morgenstern-Platz 1, 1090 Vienna, Austria\\}
\email{volker.branding@univie.ac.at}
\thanks{The author gratefully acknowledges the support of the Austrian Science Fund (FWF) through the project ``Geometric Analysis of Biwave Maps'' (DOI: 10.55776/P34853)
}

\begin{document}

\begin{abstract}
We derive various classification results for polyharmonic helices,
which are polyharmonic curves whose geodesic curvatures are all constant, in space forms.
We obtain a complete classification of triharmonic helices in
spheres of arbitrary dimension. Moreover, 
we show that polyharmonic helices of arbitrary order
with non-zero geodesic curvatures
to space forms of negative curvature must be geodesics.
\end{abstract} 

\maketitle

\section{Introduction and Results}
One of the most important objects on a given Riemannian manifold are 
\emph{geodesics} as they are the curves with minimal
distance between two given points.
The most practical way to deal with geodesics is to employ a variational approach.
Hence, we consider a curve \(\gamma\colon I\to M\),
where \(I\subset\R\) represents an interval, \((M,g)\) a Riemannian manifold
and by \(s\) we represent the parameter of the curve \(\gamma\).
Moreover, we use the notation \(\gamma'=\frac{d\gamma}{ds}\).
Then, we define the energy of the curve \(\gamma\) by
\begin{align}
\label{energy-curve}
E_1(\gamma)=E(\gamma):=\frac{1}{2}\int_I|\gamma'|^2ds.
\end{align}

The critical points of \eqref{energy-curve} are precisely \emph{geodesics}
and are characterized by the equation
\begin{align*}
\tau(\gamma):=\nabla_{\gamma'}\gamma'=0,
\end{align*}
which is a second order non-linear ordinary differential equation for the curve
\(\gamma\). The existence of geodesics on a Riemannian manifold is guaranteed
due to the celebrated Theorem of \emph{Hopf-Rinow}.

Another interesting class of curves can be obtained by 
extremizing the \emph{bienergy} of a curve \(\gamma\)
which is given by 
\begin{align}
\label{bienergy-curve}
E_{2}(\gamma):=\frac{1}{2}\int_I|\tau(\gamma)|^2ds.
\end{align}
The critical points of this energy are called \emph{biharmonic curves}
and are characterized by a non-linear ordinary differential equation
of fourth order. An important aspect of biharmonic curves is the fact
that they always have constant geodesic curvature, we refer to \cite{MR4733870} for
references and the current status of research on biharmonic curves.

A higher order generalization of both \eqref{energy-curve} and \eqref{bienergy-curve}
is provided by the \(r\)-energy of a curve
\begin{align}
\label{energy-poly-curve}
E_r(\gamma):=\frac{1}{2}\int_I|\nabla_T^{r-1}T|^2ds,
\end{align}
where \(T=\gamma'\) is the tangent vector of \(\gamma\).

The critical points of \eqref{energy-poly-curve} were calculated in \cite{MR3007953,wang} and are
characterized by the equation
\begin{align}
\label{el-poly}
0=\tau_r(\gamma)=\nabla^{2r-1}_TT+\sum_{l=0}^{r-2}(-1)^lR^M(\nabla_T^{2r-3-l}T,\nabla_T^lT)T
\end{align}
with \(R^M\) being the curvature tensor of the manifold \(M\).

Solutions of \eqref{el-poly} are called \emph{polyharmonic curves of order \(r\)} or shortly \emph{r-harmonic curves}.
In the case of \(r=1\) the energy \eqref{energy-poly-curve} reduces to the usual energy of a curve \eqref{energy-curve}
whose critical points are \emph{geodesics}. Clearly, every geodesic is a solution of the equation for \(r\)-harmonic
curves \eqref{el-poly}, hence we are interested in finding non-geodesic solutions of \eqref{el-poly}
which we will call \emph{proper} \(r\)-harmonic curves.

For the current status of research on higher order variational problems
we refer to \cite{MR4106647}, a collection of recent results 
on \(r\)-harmonic curves can be found in \cite{MR4542687}.
We want to point out that \(r\)-harmonic curves have important applications
in template matching and in the modelling of robot motions,
see for example \cite{MR2864799}.

Throughout this article we will use the following notation.
By \(s\) we will denote the parameter of the curve \(\gamma\),
the first, second and third derivative of \(\gamma\) will be written as \(T:=\gamma',\gamma''\) and \(\gamma'''\),
respectively. The \(l\)-th derivative of \(\gamma\) with respect to \(s\) will be denoted by \(\gamma^{(l)}\)
where \(l=4,\ldots,2r\). 

Moreover, we use the terminology \emph{helix} to represent
a curve whose geodesic curvatures \(k_i,i=1,\ldots\) are all constant.

Our first result is the explicit form of the Euler-Lagrange equation
for \(4\)-harmonic curves on the Euclidean sphere with the round metric.

\begin{Satz}
\label{thm:el-4-harmonic}
Let \(\gamma\colon I\to\s^n\subset\R^{n+1}\) be a curve which is parametrized by arclength.
Then \(\gamma\) is a proper \(4\)-harmonic curve if it is a non-geodesic solution of 
\begin{align}
\label{ode-4harmonic-sphere}
0=&\gamma^{(8)}
+2\gamma^{(6)}+3\gamma^{(4)}-\gamma^{(4)}|\gamma''|^2
-6\gamma''|\gamma''|^2+4\gamma''-2\gamma''\langle\gamma^{(4)},\gamma''\rangle 
\\
\nonumber & 
+5\frac{d^2}{ds^2}\big(\gamma'\langle\gamma^{(4)},\gamma'\rangle\big)-\frac{d^4}{ds^4}\big(|\gamma''|^2\gamma)
-6\gamma'\frac{d}{ds}|\gamma''|^2-2\gamma'\frac{d}{ds}\langle\gamma^{(4)},\gamma''\rangle
-5\frac{d}{ds}\big(\gamma''\langle\gamma^{(4)},\gamma'\rangle\big) \\
\nonumber&-\gamma\big(
\langle\gamma^{(8)},\gamma\rangle
+2\langle\gamma^{(6)},\gamma\rangle+3|\gamma''|^2
-|\gamma''|^4
+6|\gamma''|^2-4+2\langle\gamma^{(4)},\gamma''\rangle 
\big) \\
\nonumber&-\gamma\bigg(
5\langle\gamma,\frac{d^2}{ds^2}\big(\gamma'\langle\gamma^{(4)},\gamma'\rangle\big)\rangle
-\langle\gamma,\frac{d^4}{ds^4}\big(|\gamma''|^2\gamma)\rangle
-5\langle\gamma,\frac{d}{ds}\big(\gamma''\langle\gamma^{(4)},\gamma'\rangle\big)\rangle
\bigg).
\end{align}
\end{Satz}

\begin{Bem}
One explicit solution of \eqref{ode-4harmonic-sphere} can be obtained as follows: 
The curve \(\gamma\colon I\to\s^n\) given by
\begin{align*}
\gamma(s)=\cos(\sqrt{4}s)e_1+\sin(\sqrt{4}s)e_2+e_3,
\end{align*}
where \(e_i,i=1,2,3\) are mutually perpendicular and satisfy \(|e_1|^2=|e_2|^2=\frac{1}{4},|e_3|^2=\frac{3}{4}\),
is a proper $4$-harmonic curve which is parametrized by arclength.
The existence of this particular $4$-harmonic curve was already established 
in \cite[Theorem 1.5]{MR4542687} without using the Euler-Lagrange equation
\eqref{ode-4harmonic-sphere}.
\end{Bem}

Our next result provides a characterization of triharmonic helices in the sphere
extending the analysis presented in \cite{MR4542687}.
\begin{Satz}
\label{thm:explicit-tri}
Consider a curve \(\gamma\colon I\to\s^n\subset\R^{n+1}\) of the form
  \begin{align*}
   \gamma(s)=\cos(as)e_1+\sin(as)e_2+\cos(bs)e_3+\sin(bs)e_4
  \end{align*}
with \(|e_1|^2=|e_2|^2, |e_3|^2=|e_4|^2\) and \(|e_1|^2+|e_3|^2=1\).
Then, \(\gamma\) is a proper triharmonic curve parametrized by arclength if the following algebraic relations hold
\begin{align}
\label{eq:ab-triharmonic}
a^4+b^4-4(a^2+b^2)+3a^2b^2+3=&0,\\
\nonumber|e_1|^2a^2+|e_3|^2b^2=&1.
\end{align}
\end{Satz}

\begin{Bem}
\begin{enumerate}
\item The condition \eqref{eq:ab-triharmonic} has already been derived in \cite[Equation 2.3]{MR4542687}
using a different approach as utilized in this manuscript.
\item Setting \(a^2=x\) and \(b^2=y\) the equation \eqref{eq:ab-triharmonic} describes a particular
conic section which turns out to be a hyperbola. Hence, it is obvious that there is a whole
family of triharmonic helices in the sphere.
\end{enumerate}
\end{Bem}

It is well-known that a polyharmonic curve of order \(r\) has \(2r-2\)
non-zero geodesic curvatures and effectively lies on a target of
dimension \(2r-1\). 
Exploiting this fact the next Theorem gives some further
characterizations of \(r\)-harmonic helices in the cases \(r=3,4\).
From a computational point of view it turns out to be more effective
to work with the geodesic curvatures of a curve instead of trying
to explicitly solve the Euler-Lagrange equation.

\begin{Satz}
\label{thm:tri-four-classification}
Let \(\gamma\colon I\to M\) be a proper \(r\)-harmonic curve parametrized by arclength
where \(M\) is a space form of constant curvature \(K\).
Moreover, assume that the geodesic curvatures \(k_i,i=1,\ldots,2r-2\) are all constant.
\begin{enumerate}
 \item If \(r=3\) the geodesic curvatures \(k_j,j=1,\ldots 4\) satisfy 
 \begin{align}
\label{eq:triharmonic-curvatures}
(k_1^2+k_2^2)^2+k_2^2k_3^2&=K(2k_1^2+k_2^2),\\
\nonumber k_2k_3\big(k_1^2+k_2^2+k_3^2+k_4^2\big)&=k_2k_3K.
\end{align}
 \item If \(r=4\) the geodesic curvatures \(k_j,j=1,\ldots 6\) satisfy
 \begin{align}
 \label{eq:fourharmonic-curvatures}
(k_1^2+k_2^2)^3+&k_2^2k_3^2(2k_1^2+2k_2^2+k_3^2+k_4^2) \\
\nonumber&=K\big((k_1^2+k_2^2)^2+k_2^2k_3^2\big)+2Kk_1^2(k_1^2+k_2^2),\\
\nonumber k_2k_3\big((k_1^2+k_2^2)^2&+(k_3^2+k_4^2)^2+k_1^2k_3^2+2k_2^2k_3^2+k_4^2(k_1^2+k_2^2+k_5^2)\big) \\
\nonumber&=k_2k_3K(2k_1^2+k_2^2+k_3^2+k_4^2),\\
&\nonumber k_2k_3k_4k_5\big(k_1^2+k_2^2+k_3^2+k_4^2+k_5^2+k_6^2\big)=k_2k_3k_4k_5K.
\end{align}
\end{enumerate}
\end{Satz}

\begin{Bem}
\begin{enumerate}
\item Note that in \eqref{eq:triharmonic-curvatures} and \eqref{eq:fourharmonic-curvatures}
there does not appear a factor of \(k_1\). As we are considering proper tri- and \(4\)-harmonic curves we have that \(k_1\neq 0\) such that this factor can be split off.
\item In the case of \(k_1,k_2\neq 0\) and \(k_j=0,j\geq 3\) the first equation 
of \eqref{eq:triharmonic-curvatures} and the first equation of \eqref{eq:fourharmonic-curvatures} reduce to the formula obtained in \cite[Theorem 1.1]{MR4542687}.
\item An immediate consequence of Theorem \ref{thm:tri-four-classification} is that
triharmonic and \(4\)-harmonic helices whose geodesic curvatures are all non-zero need to be geodesics if the target is a space form of non-positive curvature. 
More precisely, for triharmonic helices the first equality of \eqref{eq:triharmonic-curvatures} leads to a contradiction if \(K\) is non-positive,
while for \(4\)-harmonic helices \eqref{eq:fourharmonic-curvatures} gives a contradiction whenenver \(K\) is non-positive.
\end{enumerate}
\end{Bem}

Employing Theorem \ref{thm:tri-four-classification} we can deduce the following
classification result:
\begin{Satz}
\label{thm:triharmonic-helices}
A proper triharmonic helix \(\gamma\colon I\to\s^n\subset\R^{n+1}\)
parametrized by arclength
must be one of the following:
\begin{enumerate}
\item A circle of the form ($n\geq 2$)
\begin{align*}
\gamma(s)=\cos(\sqrt{3}s)e_1+\sin(\sqrt{3}s)e_2+e_3,
\end{align*}
where \(e_i,i=1,2,3\) are mutually perpendicular and satisfy \(|e_1|^2=|e_2|^2=\frac{1}{3},|e_3|^2=\frac{2}{3}\).

\item A non-planar curve of the form ($n\geq 3$)
\begin{align*}  
  \gamma(s)=\cos(as)e_1+\sin(as)e_2+\cos(bs)e_3+\sin(bs)e_4
  \end{align*}
with \(|e_1|^2=|e_2|^2, |e_3|^2=|e_4|^2, |e_1|^2+|e_3|^2=1\)
and \(e_i,i=1,\ldots,4\) mutually perpendicular
satisfying
\begin{align*}
a^4+b^4-4(a^2+b^2)+3a^2b^2+3=&0, \qquad |e_1|^2a^2+|e_3|^2b^2=1.
\end{align*}
\end{enumerate}
\end{Satz}

\begin{Bem}
It is quite remarkable that triharmonic helices on the sphere
have the same structure as biharmonic curves on the sphere.
In both cases there exist two families of the form detailed 
in the previous theorem, which consist of a planar family
and a non-planar generalization, see Theorem \ref{prop-cmo}
for the precise details on biharmonic curves.
However, biharmonic curves necessarily have constant curvature
while there may be triharmonic curves on the sphere of non-constant
geodesic curvature. In order to obtain a complete classification
of triharmonic curves on the sphere
one would need to obtain a full understanding of the non-constant
curvature case as well.
\end{Bem}

The last result of this manuscript provides a characterization of 
polyharmonic helices of arbitrary order
whose geodesic curvatures are all different from zero.

\begin{Satz}
\label{thm:structure}
Let \(\gamma\colon I\to M\) be an \(r\)-harmonic curve parametrized by arclength
whose geodesic curvatures \(k_i,i=1,\ldots,2r-2\) are all constant and non-zero.
Moreover, suppose that \(M\) is a space form
of constant curvature \(K\). Then, the following equation holds
\begin{align*}
\sum_{j=1}^{2r-2}k_j^2=K.
\end{align*}
\end{Satz}

\begin{Bem}
The previous Theorem gives further insights into the structure of higher order variational problems. 
In particular, it gives further evidence to support the claim that polyharmonic maps
to space forms of negative curvature must be harmonic,
by showing that there do not exist proper \(r\)-harmonic helices when \(K\) is non-positive,
while there may be additional solutions in the case of a spherical target. 
So far, this observations was mostly made in the case of codimension one, that is for polyharmonic hypersurfaces, see for example \cite{MR4462636}.
The analysis above suggests that this fact stays true in higher codimension.
\end{Bem}

Furthermore, we will collect a number of results giving rise to the following
conjecture:

\begin{Con}
\label{conjecture}
The equation for \(r\)-harmonic curves \eqref{el-poly} admits
solutions with non-constant geodesic curvature 
\begin{align}
k_1(s)=\frac{\alpha}{s^{r-2}},\qquad \alpha\in\R,\alpha\neq 0,
\end{align}
where \(s\) represents the parameter of the curve \(\gamma\)
which we assume to be parametrized by arclength.
\end{Con}

This conjecture is based on
\begin{enumerate}
 \item the well-known fact that biharmonic curves
(\(r=2\)) necessarily have constant geodesic curvature,
\item the results of the recent article on triharmonic curves \cite{MR4308322},
\item and the observations presented in Subsection \ref{subsection-conjecture}.
\end{enumerate}
\par\medskip
This article is organized as follows:
In Section 2 we provide some background material on the Euler-Lagrange method
and use it to reprove a number of well-known results on biharmonic curves
on the Euclidean sphere. Finally, in Section 3 we give the proofs of the main results of the article.

\section{Some preliminary results}
Throughout this article we consider a space form of constant curvature 
\(K\) in which case the Riemann curvature tensor acquires the simple form
\begin{align*}
R(X,Y)Z=K(\langle Y,Z\rangle X-\langle X,Z\rangle Y),
\end{align*}
where \(X,Y,Z\) are vector fields and \(K\)
represents the constant curvature of the space form.

In this case the equation for polyharmonic curves \eqref{el-poly} simplifies to
\begin{align}
\label{polyharmonic-spaceform}
\tau_r(\gamma)=\nabla_T^{2r-1}T+K\sum_{l=0}^{r-2}(-1)^l\big(\langle T,\nabla^l_TT\rangle\nabla_T^{2r-3-l}T
-\langle T,\nabla_T^{2r-3-l}T\rangle\nabla_T^lT\big).
\end{align}

\subsection{The Euler-Lagrange method for polyharmonic curves}
Let us briefly recall the so-called \emph{Euler-Lagrange method}
which is a powerful tool in the analysis of one-dimensional variational problems.
This method is the cornerstone of the Lagrangian formulation of classical mechanics
in theoretical physics, see for example \cite[Chapter 7]{MR1723696}.
Moreover, this method can also successfully be applied in order to study biharmonic curves \cite{MR2274737},
biharmonic maps \cite{MR3577677,MR3045700} and also higher order variational problems \cite[Theorem 4.5]{MR4106647}.

The following theorem may be well-known in the mathematics community.
However, for the sake of completeness we also provide a  
complete proof below.

\begin{Satz}
\label{theorem-el-poly}
Let \(\gamma\colon I\to\R^q\) be a curve.
Suppose we have an energy functional
\begin{align*}
E_r(\gamma)=\int_I\mathcal{L}_r~ds,
\end{align*}
where the Lagrangian 
\begin{align*}
\mathcal{L}_r=\mathcal{L}_r(\gamma,\gamma',\ldots,\gamma^{(r-1)},\gamma^{(r)})
\end{align*}
may depend on the derivatives of the curve \(\gamma\) up to order \(r\).

Then, \(\gamma\) is a critical point of \(E_r(\gamma)\) if the following
ordinary differential equation holds
\begin{align}
\sum_{l=1}^r(-1)^l\frac{d^l}{ds^l}\frac{\partial\mathcal{L}_r}{\partial\gamma^{(l)}}+\frac{\partial\mathcal{L}_r}{\partial\gamma}=0.
\end{align}
\end{Satz}

\begin{proof}
We choose \(\beta\in C^\infty_c(I,\R^q)\) and compute the first variation of \(E_r(\gamma)\) as follows
\begin{align*}
\frac{d}{dt}\big|_{t=0}&E_r(\gamma+t\beta)\\
&=\frac{d}{dt}\int_I\mathcal{L}_r(\gamma+t\beta,\gamma'+t\beta',\gamma''+t\beta'',\ldots,\gamma^{(r-1)}+t\beta^{(r-1)},\gamma^{(r)}+t\beta^{(r)})~ds\big|_{t=0} \\
&=\int_I\frac{d}{dt}\big(\mathcal{L}_r(\gamma+t\beta,\gamma'+t\beta',\gamma''+t\beta'',\ldots,\gamma^{(r-1)}+t\beta^{(r-1)},\gamma^{(r)}+t\beta^{(r)})\big)~ds\big|_{t=0} \\
&=\int_I\big(\frac{\partial\mathcal{L}_r}{\partial\gamma}\beta+\frac{\partial\mathcal{L}_r}{\partial\gamma'}\beta'
+\frac{\partial\mathcal{L}_r}{\partial\gamma''}\beta''+\ldots+\frac{\partial\mathcal{L}_r}{\partial\gamma^{(r-1)}}\beta^{(r-1)}
+\frac{\partial\mathcal{L}_r}{\partial\gamma^{(r)}}\beta^{(r)}\big)~ds.
\end{align*}
Now, we use integration by parts
\begin{align*}
\int_I\frac{\partial\mathcal{L}_r}{\partial\gamma^{(p)}}\beta^{(p)}~ds
=\int_I(-1)^p\frac{d^p}{ds^p}\big(\frac{\partial\mathcal{L}_r}{\partial\gamma^{(p)}}\big)\beta~ ds,\qquad 1\leq p\leq r,
\end{align*}
where we used that \(\beta\) is compactly supported. 
By combining both equations the proof is complete.
\end{proof}

In the following we will often make use of the following lemma
which follows from a direct calculation.
\begin{Lem}
\label{lemma-identities-sphere}
Let \(\gamma\colon I\to\s^n\subset{\R}^{n+1}\) be a curve which is parametrized by arclength.
Then the following identities hold
\begin{align*}
\langle\gamma,\gamma'\rangle&=0,\qquad \langle\gamma'',\gamma\rangle=-1,\qquad \langle\gamma''',\gamma\rangle=0,\qquad
\langle\gamma',\gamma''\rangle=0,\\
\langle\gamma^{(4)},\gamma\rangle+\langle\gamma''',\gamma'\rangle&=0,\qquad
\langle\gamma^{(4)},\gamma\rangle=|\gamma''|^2.
\end{align*}
\end{Lem}

Throughout this section we will frequently make use of the inclusion map \(\iota\colon\s^n\to\R^{n+1}\)
and also exploit the special structure of the Levi-Civita connection on the sphere
\begin{align*}
d\iota(\nabla_TX)=X'+\langle X,\gamma'\rangle\gamma,
\end{align*}
where \(X\) is a vector field on \(\s^n\subset\R^{n+1}\).

\subsection{Biharmonic curves on the sphere}
In order to highlight the power of the Euler-Lagrange method we will first investigate biharmonic curves on the Euclidean sphere and give a new-proof of some well-known results
which serves as an inspiration for the classification results on triharmonic helices
presented in this manuscript.

The intrinsic form of the equation for biharmonic curves on the sphere is given by
\begin{align*}
\tau_2(\gamma)=\nabla^3_TT+|T|^2\nabla_TT-\langle T,\nabla_TT\rangle T.
\end{align*}

Assuming that \(\s^n\subset\R^{n+1}\) and considering a curve \(\gamma\colon I\to\s^n\subset\R^{n+1}\)
we obtain the following Lagrangian for biharmonic curves on the sphere
\begin{align}
\label{lagrangian-biharmonic}
\mathcal{L}^{\s^n}_2(\gamma'',\gamma',\gamma)=|\gamma''|^2-|\gamma'|^4+\lambda (|\gamma|^2-1).
\end{align}

Note that we have to include the Lagrange multiplyer \(\lambda\) as we are constraining
the curve \(\gamma\) to be on the unit sphere.

Then, employing Theorem \ref{theorem-el-poly}, a direct calculation shows that
the critical points of \eqref{lagrangian-biharmonic} are given by
\begin{align*}
\frac{d^2}{ds^2}\big(\frac{\partial\mathcal{L}^{\s^n}_{2}}{\partial\gamma''}\big)
-\frac{d}{ds}\big(\frac{\partial\mathcal{L}^{\s^n}_{2}}{\partial\gamma'}\big)
+\frac{\partial\mathcal{L}^{\s^n}_{2}}{\partial\gamma}
=2(\gamma^{(4)}
+2(|\gamma'|^2)'\gamma'+2|\gamma'|^2\gamma''+\lambda\gamma).
\end{align*}

Taking also into account the variation of \(\mathcal{L}^{\s^n}_2\) with respect to the Lagrange multiplyer \(\lambda\)
we obtain the Euler-Lagrange equation
\begin{align}
\label{el-sphere-biharmonic-a}
\gamma^{(4)}+2(|\gamma'|^2)'\gamma'+2|\gamma'|^2\gamma''+\lambda\gamma=0
\end{align}
together with the constraint \(|\gamma|^2=1\).

From now on, we will assume that the curve \(\gamma\) is parametrized with respect to arclength,
that is \(|\gamma'|^2=1\), such that \eqref{el-sphere-biharmonic-a} simplifies to
\begin{align}
\label{el-sphere-biharmonic-b}
\gamma^{(4)}+2\gamma''+\lambda\gamma=0.
\end{align}

In order to determine \(\lambda\) we test \eqref{el-sphere-biharmonic-b} with \(\gamma\)
and find
\begin{align*}
\lambda=-\langle\gamma^{(4)},\gamma\rangle-2\langle\gamma,\gamma''\rangle=-|\gamma''|^2+2,
\end{align*}
where we used the identities provided by Lemma \ref{lemma-identities-sphere}
and thus exploited the fact that \(\gamma\) is parametrized with respect to arclength.

It is well-known that biharmonic curves have constant geodesic curvature
which, using our framework, can be seen as follows:
\begin{Bem}
\begin{enumerate}
\item It is easy to see that the Lagrange multiplyer \(\lambda\) and the geodesic curvature \(k_1\) of the curve
\(\gamma\) are related via the identity
\begin{align*}
\lambda=-|\gamma''|^2+2=-k_1^2+1.
\end{align*}
Hence, the inclusion of the Lagrange multiplyer \(\lambda\) in \eqref{lagrangian-biharmonic} has the effect that
it forces the curve \(\gamma\) to have constant geodesic curvature.

This fact is well-known and is usually deduced by choosing a Frenet-frame for the curve \(\gamma\)
and analyzing the associated Frenet equations.

 \item We will present another short argument why biharmonic curves always need to have constant geodesic
curvature that holds for biharmonic curves on an arbitrary manifold.

Suppose we have a biharmonic curve on a Riemannian manifold, then it satisfies
\begin{align*}
\tau_2(\gamma)=\nabla^3_TT+R(\nabla_TT,T)T=0.
\end{align*}
Multiplying this equation with \(T\) we obtain
\begin{align*}
\langle\nabla^3_TT,T\rangle=0
\end{align*}
which implies
\begin{align*}
\frac{d}{ds}\langle\nabla^2_TT,T\rangle-\frac{1}{2}\frac{d}{ds}|\nabla_TT|^2=0.
\end{align*}
Exploiting that the curve is parametrized with respect
to arclength we can then deduce that 
\begin{align*}
\frac{d}{ds}|\nabla_TT|^2=0,
\end{align*}
which implies that \(k_1^2=const\).
\end{enumerate}
\end{Bem}

Combining the previous observations we get the following well-known result:
\begin{Prop}
Let \(\gamma\colon I\to\s^n\subset\R^{n+1}\) be a curve parametrized by arclength.
Then \(\gamma\) is biharmonic if
\begin{align}
\label{biharmonic-sphere-ode}
\gamma^{(4)}+2\gamma''+\gamma(2-|\gamma''|^2)=0.
\end{align}
\end{Prop}

\begin{Bem}
\begin{enumerate}
 \item The equation for biharmonic curves on spheres \eqref{biharmonic-sphere-ode} 
 was first derived in \cite[Corollary 4.2]{MR1919374} making
use of geometric methods. In that reference the equation for biharmonic curves to spheres
is given in the following form
\begin{align*}
\gamma^{(4)}+2\gamma''+\gamma(1-k_1^2)=0,
\end{align*}
where \(k_1\) represents the geodesic curvature of the curve \(\gamma\).
Noting that
\begin{align*}
k_1^2=|\nabla_{T}T|^2=|\gamma''|^2-|\gamma'|^4=|\gamma''|^2-1
\end{align*}
it is obvious that this version is the same as \eqref{biharmonic-sphere-ode}.
\item We would like to point out that it is necessary to include the Lagrange multiplyer \(\lambda\)
in the Lagrangian \eqref{lagrangian-biharmonic} as we are dealing with a constraint variational problem.
However, we do not have to include a second Lagrange multiplyer to justify
that the curve is parametrized with respect to arclength. The fact that we are choosing
an arclength parametrization can be considered as making a convenient choice 
in order to simplify our calculations but it is not a constraint required from
the actual variational problem.
\end{enumerate}

\end{Bem}

The following result was proved in 
\cite{MR1863283} and 
\cite[Proposition 4.4]{MR1919374}.

\begin{Satz}
\label{prop-cmo}
Let \(\gamma\colon I\to\s^n\subset\R^{n+1}\) be a curve parametrized by arclength.
Then there exist the following two classes of proper biharmonic curves on \(\s^n\):
\begin{enumerate}
 \item When \(k_1^2=1\) these are circles parametrized by
\begin{align}
\label{biharmonic-sphere-a}
\gamma(s)=\cos(\sqrt{2}s)e_1+\sin(\sqrt{2}s)e_2+e_3,
\end{align}
where \(e_i,i=1,2,3\) are constant orthogonal vectors satisfying 
\(|e_1|^2=|e_2|^2=|e_3|^2=\frac{1}{2}\).
\item When \(0<k_1^2<1\) they are non-planar curves parametrized as follows
\begin{align}
\label{biharmonic-sphere-b}
\gamma(s)=\cos(as)e_1+\sin(as)e_2
+\cos(bs)e_3+\sin(bs)e_4,
\end{align}
where \(|e_i|^2=\frac{1}{2},i=1,\ldots,4\) and \(a^2+b^2=2\) with \(a^2\neq b^2\).
\end{enumerate}
\end{Satz}

In the following we will give a different proof of Theorem \ref{prop-cmo} as was 
originally presented in \cite{MR1863283,MR1919374}, making use of the Euler-Lagrange method, as we want to employ it frequently in the rest of this article.

\begin{proof}[Proof of Theorem \ref{prop-cmo}]
In order to find the first class of solutions \eqref{biharmonic-sphere-a} we make the ansatz
\begin{align*}
\gamma(s)=\cos(as)e_1+\sin(as)e_2+e_3,
\end{align*}
where \(e_i,i=1,2,3\) are constant orthogonal vectors satisfying \(|e_1|^2=|e_2|^2\)
and \(|e_1|^2+|e_3|^2=1\) as we require \(|\gamma|^2=1\) and \(a\in\R\).
In the following we set \(\alpha^2:=|e_1|^2.\)
Inserting this ansatz into the Lagrangian for biharmonic curves \eqref{lagrangian-biharmonic} we find
\begin{align*}
\mathcal{L}^{\s^n}_2(\alpha)=a^4(\alpha^2-\alpha^4).
\end{align*}
To determine the critical points of \(\mathcal{L}^{\s^n}_2(\alpha)\) we calculate
\begin{align*}
\frac{d}{d\alpha}\mathcal{L}^{\s^n}_2(\alpha)=2a^4\alpha(1-2\alpha^2)
\end{align*}
and it is clear that this expression vanishes if \(\alpha^2=\frac{1}{2}\).
Finally, we use the fact that \(\gamma\) is parametrized with respect to arclength
which, given our ansatz, is expressed via \(\alpha^2a^2=1\).
Hence, we get that \(a^2=2\) completing the proof.

In order to obtain the second class of solutions \eqref{biharmonic-sphere-b} we make the ansatz
\begin{align*}
\gamma(s)=\cos(as)e_1+\sin(as)e_2+\cos(bs)e_3+\sin(bs)e_4,
\end{align*}
where \(|e_1|^2=|e_2|^2,|e_3|^2=|e_4|^2\) .
We set \(\alpha_j^2:=|e_j|^2, j=1,3.\)
Inserting this ansatz into the Lagrangian for biharmonic curves \eqref{lagrangian-biharmonic} we get
\begin{align*}
\mathcal{L}^{\s^n}_2(\alpha_1,\alpha_3,\lambda)=
a^4(\alpha_1^2-\alpha_1^4)+b^4(\alpha_3^2-\alpha_3^4)
-2a^2b^2\alpha_1^2\alpha_3^2+\lambda(\alpha_1^2+\alpha_3^2-1).
\end{align*}
In order to find the critical points of this Lagrangian we differentiate with respect to \(\alpha_1,\alpha_3,\lambda\)
and set the resulting equations equal to zero leading to the system
\begin{align*}
a^4(1-2\alpha_1^2)-2a^2b^2\alpha_3^2-\lambda=&0,\\
b^4(1-2\alpha_3^2)-2a^2b^2\alpha_1^2-\lambda=&0,\\
\alpha_1^2+\alpha_3^2=&1.
\end{align*}
Combining this set of equations we find after some algebraic manipulations
\begin{align*}
(a^4+b^4-2a^2b^2)(1-2\alpha_1^2)=0
\end{align*}
from which we deduce that \(\alpha_1^2=\alpha_3^2=\frac{1}{2}\).
The requirement that \(\gamma\) is parametrized with respect to arclength is then given by the constraint
\(a^2+b^2=2,a^2\neq b^2\) completing the proof.
\end{proof}

\section{Proofs of the main results}
In this section we provide the proofs of the main results of this article.
\begin{proof}[Proof of Theorem \ref{thm:el-4-harmonic}]
The proof is based on the Euler-Lagrange method provided by Theorem \ref{theorem-el-poly}.
Recall that the 4-energy of a curve \(\gamma\colon I\to\s^n\) is given by
\begin{align*}
E_4(\gamma)=\int_I|\nabla^3_TT|^2ds.
\end{align*}
We again make use of the embedding \(\iota\colon\s^n\to\R^{n+1}\) which helps us
rewrite
\begin{align*}
d\iota(\nabla^3_TT)=\gamma^{(4)}+4\langle\gamma''',\gamma'\rangle\gamma+3|\gamma''|^2\gamma
+5\langle\gamma',\gamma''\rangle\gamma'+|\gamma'|^2\gamma''+|\gamma'|^4\gamma.
\end{align*} 
Consequently, the Lagrangian associated with the \(4\)-energy for a curve
on the sphere has the form
\begin{align}
\label{lagrangian-fourharmonic}
\mathcal{L}^{\s^n}_4(\gamma^{(4)},\gamma''',\gamma'',\gamma',\gamma)=&
|\gamma^{(4)}|^2+16|\langle\gamma''',\gamma'\rangle|^2+9|\gamma''|^4
+35|\langle\gamma'',\gamma'\rangle|^2|\gamma'|^2 \\
\nonumber&+|\gamma'|^4|\gamma''|^2-|\gamma'|^8 \\
\nonumber&+8\langle\gamma''',\gamma'\rangle\langle\gamma^{(4)},\gamma\rangle
+6\langle\gamma^{(4)},\gamma\rangle|\gamma''|^2
+10\langle\gamma'',\gamma'\rangle\langle\gamma^{(4)},\gamma'\rangle \\
\nonumber&+2|\gamma'|^2\langle\gamma^{(4)},\gamma''\rangle+2|\gamma'|^4\langle\gamma^{(4)},\gamma\rangle \\
\nonumber&+24\langle\gamma''',\gamma'\rangle|\gamma''|^2 \\
\nonumber&+\lambda(|\gamma|^2-1),
\end{align}
where we again introduced the Lagrange multiplyer \(\lambda\) to constrain the curve \(\gamma\) to \(\s^n\).

By a direct calculation we find taking into account that the curve \(\gamma\)
is parametrized with respect to arclength
\begin{align*}
\frac{d^4}{ds^4}\frac{\partial\mathcal{L}^{\s^n}_4}{\partial\gamma^{(4)}}=&2\gamma^{(8)}
-2\frac{d^4}{ds^4}\big(|\gamma''|^2\gamma)+2\gamma^{(6)}+2\gamma^{(4)}
, \\
\frac{d^3}{ds^3}\big(\frac{\partial\mathcal{L}^{\s^n}_{4}}{\partial\gamma'''}\big)
=& 0,\\
\frac{d^2}{ds^2}\big(\frac{\partial\mathcal{L}^{\s^n}_{4}}{\partial\gamma''}\big)
=& 2\gamma^{(6)}+2\gamma^{(4)}+10\frac{d^2}{ds^2}\big(\gamma'\langle\gamma^{(4)},\gamma'\rangle\big),\\
\frac{d}{ds}\big(\frac{\partial\mathcal{L}^{\s^n}_{4}}{\partial\gamma'}\big)
=&12\gamma''|\gamma''|^2-8\gamma''+4\gamma''\langle\gamma^{(4)},\gamma''\rangle \\
&+12\gamma'\frac{d}{ds}|\gamma''|^2+4\gamma'\frac{d}{ds}\langle\gamma^{(4)},\gamma''\rangle
+10\frac{d}{ds}\big(\gamma''\langle\gamma^{(4)},\gamma'\rangle\big)
,\\
\frac{\partial\mathcal{L}^{\s^n}_{4}}{\partial\gamma}
=& -2\gamma^{(4)}|\gamma''|^2+2\gamma^{(4)}+2\lambda\gamma.
\end{align*}
Varying \eqref{lagrangian-fourharmonic} with respect to the Lagrange multiplyer \(\lambda\) we obtain
the constraint \(|\gamma|^2=1\).
Hence, from Theorem \ref{theorem-el-poly} we can deduce that
\begin{align*}
0=&\gamma^{(8)}
+2\gamma^{(6)}+3\gamma^{(4)}-\gamma^{(4)}|\gamma''|^2
-6\gamma''|\gamma''|^2+4\gamma''-2\gamma''\langle\gamma^{(4)},\gamma''\rangle 
\\
& 
+5\frac{d^2}{ds^2}\big(\gamma'\langle\gamma^{(4)},\gamma'\rangle\big)-\frac{d^4}{ds^4}\big(|\gamma''|^2\gamma)
-6\gamma'\frac{d}{ds}|\gamma''|^2-2\gamma'\frac{d}{ds}\langle\gamma^{(4)},\gamma''\rangle
-5\frac{d}{ds}\big(\gamma''\langle\gamma^{(4)},\gamma'\rangle\big) 
+\lambda\gamma.
\end{align*}
In order to determine \(\lambda\) we form the scalar product with \(\gamma\),
using the identifies provided by Lemma \ref{lemma-identities-sphere}
and inserting back 
into the above equation completes the proof.
\end{proof}

As the proof of Theorem \ref{thm:explicit-tri} is based on the 
Lagrangian for triharmonic curves 
we again use the embedding of \(\s^n\) into \(\R^{n+1}\) via the map \(\iota\) and find
\begin{align*}
d\iota(\nabla^2_TT)=\gamma'''+3\langle\gamma'',\gamma'\rangle\gamma+|\gamma'|^2\gamma'.
\end{align*}

Thus, we obtain the following Lagrangian
\begin{align}
\label{lagrangian-triharmonic}
\mathcal{L}_3^{\s^n}(\gamma,\gamma',\gamma'',\gamma''')=|\gamma'''|^2+9|\langle\gamma'',\gamma'\rangle|^2+|\gamma'|^6
+6\langle\gamma'',\gamma'\rangle\langle\gamma''',\gamma\rangle+2|\gamma'|^2\langle\gamma',\gamma'''\rangle
+\lambda(|\gamma|^2-1)
\end{align}
with the Lagrange multiplyer \(\lambda\in\R\).

In order to prove Theorem \ref{thm:explicit-tri} we first establish the following 
\begin{Prop}
Consider a curve \(\gamma\colon I\to\s^n\) of the form
  \begin{align}
  \label{ansatz-triharmonic-d}
   \gamma(s)=\cos(as)e_1+\sin(as)e_2+\cos(bs)e_3+\sin(bs)e_4
  \end{align}
with \(|e_1|^2=|e_2|^2, |e_3|^2=|e_4|^2\).
Then, \(\gamma\) is a proper triharmonic curve parametrized by arclength if the following algebraic relations hold
\begin{align}
\label{algebraic-system-triharmonic}
a^6(1-2\alpha_1^2)-2a^4+3a^2-2a^2b^4\alpha_3^2+\lambda=&0,\\
\nonumber b^6(1-2\alpha_3^2)-2b^4+3b^2-2b^2a^4\alpha_1^2+\lambda=&0,\\
\nonumber a^2\alpha_1^2+b^2\alpha_3^2=&1,\\
\nonumber \alpha_1^2+\alpha_3^2=&1,
\end{align}
whenever \(a^2,b^2\neq 1\).
Here, \(\lambda\in\R\) and we have set \(\alpha_j^2=|e_j|^2,j=1,3\).
\end{Prop}
\begin{proof}
From the ansatz \eqref{ansatz-triharmonic-d} we get 
\begin{align*}
|\gamma^{(l)}|^2=\alpha_1^2a^{2l}+\alpha_3^2b^{2l},\qquad l=1,2,3.
\end{align*}
Inserting into the Lagrangian for triharmonic curves \eqref{lagrangian-triharmonic} we 
obtain 
\begin{align*}
\mathcal{L}^{\s^n}_3(\alpha_1,\alpha_3,\lambda)=&\alpha_1^2a^{6}+\alpha_3^2b^{6}
+(\alpha_1^2a^{2}+\alpha_3^2b^{2})^3 
-2(\alpha_1^2a^{2}+\alpha_3^2b^{2})(\alpha_1^2a^{4}+\alpha_3^2b^{4}) \\
&+\lambda(\alpha_1^2+\alpha_3^2-1)\\
=&a^6(\alpha_1^2-2\alpha_1^4+\alpha_1^6)+b^6(\alpha_3^2-2\alpha_3^4+\alpha_3^6)\\
&+a^2b^4(3\alpha_1^2\alpha_3^4-2\alpha_1^2\alpha_3^2)
+a^4b^2(3\alpha_1^4\alpha_3^2-2\alpha_1^2\alpha_3^2)\\
&
+\lambda(\alpha_1^2+\alpha_3^2-1).
\end{align*}

The critical points of \(\mathcal{L}^{\s^n}_3(\alpha_1,\alpha_3,\lambda)\) are given by the set of equations
\begin{align}
\label{ansatz-triharmonic-d-critical}
0=&a^6(1-4\alpha_1^2+3\alpha_1^4)+a^2b^4(3\alpha_3^4-2\alpha_3^2)
+a^4b^2(6\alpha_1^2\alpha_3^2-2\alpha_3^2)
+\lambda,\\
\nonumber 0=&b^6(1-4\alpha_3^2+3\alpha_3^4)+a^2b^4(6\alpha_1^2\alpha_3^2-2\alpha_1^2)+a^4b^2(3\alpha_1^4-2\alpha_1^2)+\lambda,\\
\nonumber 1=&\alpha_1^2+\alpha_3^2.
\end{align}
In addition, we have the following constraint due to the requirement that
the curve \(\gamma\) is supposed to be parametrized by arclength
\begin{align*}
a^2\alpha_1^2+b^2\alpha_3^2=1.
\end{align*}
Using this constraint in the first equation of \eqref{ansatz-triharmonic-d-critical} we manipulate
\begin{align*}
0=&a^6(1-4\alpha_1^2+3\alpha_1^4)+a^4((6\alpha_1^2-2)(1-a^2\alpha_1^2)\big) 
+a^2\big(3(1-a^2\alpha_1^2)^2-2b^4\alpha_3^2\big)+\lambda \\
=&a^6(1-2\alpha_1^2)-2a^4+3a^2-2a^2b^4\alpha_3^2+\lambda.
\end{align*}
This shows the validity of the first equation in \eqref{algebraic-system-triharmonic},
the second one can be derived by the same method.
The last two equations in \eqref{algebraic-system-triharmonic} represent the fact
that \(|\gamma|^2=|\gamma'|^2=1\).
\end{proof}

\begin{Bem}
One can easily check that \(a=b=1\) solves the system \eqref{algebraic-system-triharmonic}
which corresponds to a geodesic solution.
\end{Bem}

\begin{proof}[Proof of Theorem \ref{thm:explicit-tri}]
Using the first two equations of \eqref{algebraic-system-triharmonic}
we obtain
\begin{align*}
a^6-b^6-2(a^4-b^4)+3(a^2-b^2)-2\alpha_1^2a^4(a^2-b^2)
-2\alpha_3^2b^4(a^2-b^2)=0.
\end{align*}
Employing the identity
\begin{align*}
a^6-b^6=(a^2-b^2)(a^4+b^4+a^2b^2)
\end{align*}
and assuming that \(a\neq b\) we can thus deduce
\begin{align*}
a^4+b^4+a^2b^2-2(a^2+b^2)+3-2\alpha_1^2a^4-2\alpha_3^2b^4=0.
\end{align*}
In order to manipulate the last two terms involving \(\alpha_1^2,\alpha_3^2\) we make
use of the last two equations of \eqref{algebraic-system-triharmonic} as follows
\begin{align*}
\alpha_1^2a^4+\alpha_3^2b^4=&a^2(1-\alpha_3^2b^2)+b^2(1-\alpha_1^2a^2)\\
=&a^2+b^2-a^2b^2(\alpha_1^2+\alpha_3^2) \\
=&a^2+b^2-a^2b^2
\end{align*}
such that we obtain
\begin{align*}
a^4+b^4+a^2b^2-2(a^2+b^2)+3-2\alpha_1^2a^4-2\alpha_3^2b^4=
a^4+b^4-4(a^2+b^2)+3a^2b^2+3
\end{align*}
yielding the claim.
\end{proof}

For the further analysis we recall the following
\begin{Dfn}[Frenet-frame]
Let \(\gamma\colon I\to M\) be a curve which is parametrized with respect to arclength.
Then, its Frenet-frame is defined by
\begin{align}
\label{frenet-frame}
F_1=&T,\\
\nonumber\nabla_{T} F_1=&k_1F_2,\\
\nonumber\nabla_{T} F_i=&-k_{i-1}F_{i-1}+k_iF_{i+1},\qquad i=2,\ldots,n-1,\\
\nonumber\nonumber\vdots \\
\nonumber\nabla_{T} F_n=&-k_{n-1}F_{n-1},
\end{align}
where \(k_i,i=1,\ldots n-1\) represent the curvatures of the curve \(\gamma\).
\end{Dfn}

\begin{proof}[Proof of Theorem \ref{thm:tri-four-classification}]
In order to prove the first result concerning the classification of triharmonic helices
in space forms we note that the equation for triharmonic curves in space forms
reads as
\begin{align}
\label{eq:triharmonic-space-form}
\nabla^5_TT+K\nabla^3_TT-K\langle T,\nabla^3_TT\rangle T
-K\langle T,\nabla_TT\rangle\nabla^2_TT+K\langle T,\nabla^2_TT\rangle\nabla_TT=0.
\end{align}
A direct calculation using \eqref{frenet-frame}, assuming that \(k_i,i=1,\ldots,4\)
are constant and \(k_i=0,i\geq 5\), 
shows that
\begin{align*}
\nabla^2_TT=&-k_1^2T+k_1k_2F_3,\\
\nabla^3_TT=&-k_1(k_1^2+k_2^2)F_2+k_1k_2k_3F_4,\\
\nabla^4_TT=&k_1^2(k_1^2+k_2^2)T-k_1k_2(k_1^2+k_2^2+k_3^2)F_3+k_1k_2k_3k_4F_5,\\
\nabla^5_TT=&k_1\big((k_1^2+k_2^2)^2+k_2^2k_3^2\big)F_2
-k_1k_2k_3(k_1^2+k_2^2+k_3^2+k_4^2)F_4.
\end{align*}
Inserting these identities into \eqref{eq:triharmonic-space-form} then yields
\begin{align*}
k_1\big((k_1^2+k_2^2)^2+k_2^2k_3^2-K(2k_1^2+k_2^2)\big)F_2
-k_1k_2k_3(k_1^2+k_2^2+k_3^2+k_4^2-K)F_4=0.
\end{align*}
Testing this equation with both \(F_2,F_4\) then completes the first part of the proof.

Concerning the second claim of the theorem, which is the classification of \(4\)-harmonic 
helices in space forms, we recall that in this case
the equation for \(4\)-harmonic
curves acquires the form
\begin{align}
\label{4-harmonic-space-form}
\nabla^7_TT+&K\nabla^5_TT-K\langle\nabla^5_TT,T\rangle T
-K\langle\nabla_TT,T\rangle\nabla^4_TT
+K\langle\nabla^4_TT,T\rangle\nabla_TT\\
\nonumber&+K\langle\nabla^2_TT,T\rangle\nabla^3_TT
-K\langle T,\nabla^3_TT\rangle\nabla^2_TT=0,
\end{align}
which is precisely \eqref{polyharmonic-spaceform} for \(r=4\).
In order to characterize the solutions of \eqref{4-harmonic-space-form}
with constant curvatures we use the Frenet equations \eqref{frenet-frame}
and a direct calculation shows
\begin{align*}
\nabla^2_TT=&-k_1^2T+k_1k_2F_3,\\
\nabla^3_TT=&-k_1(k_1^2+k_2^2)F_2+k_1k_2k_3F_4,\\
\nabla^4_TT=&k_1^2(k_1^2+k_2^2)T-k_1k_2(k_1^2+k_2^2+k_3^2)F_3+k_1k_2k_3k_4F_5,\\
\nabla^5_TT=&k_1\big((k_1^2+k_2^2)^2+k_2^2k_3^2\big)F_2
-k_1k_2k_3(k_1^2+k_2^2+k_3^2+k_4^2)F_4
+k_1k_2k_3k_4k_5F_6,\\
\nabla^6_TT=&
-k_1^2\big((k_1^2+k_2^2)^2+k_2^2k_3^2\big)T
+k_1k_2\big((k_1^2+k_2^2)^2+k_3^2(k_1^2+2k_2^2+k_3^2+k_4^2)\big)F_3 \\
&-k_1k_2k_3k_4(k_1^2+k_2^2+k_3^2+k_4^2+k_5^2)F_5
+k_1k_2k_3k_4k_5k_6F_7,\\
\nabla^7_TT=&-k_1\big((k_1^2+k_2^2)^3+k_2^2k_3^2(2k_1^2+2k_2^2+k_3^2+k_4^2)\big)F_2\\
&+k_1k_2k_3\big((k_1^2+k_2^2)^2+(k_3^2+k_4^2)^2
+k_1^2k_3^2+2k_2^2k_3^2+k_4^2(k_1^2+k_2^2+k_5^2)
\big)F_4 
\\
&-k_1k_2k_3k_4k_5(k_1^2+k_2^2+k_3^2+k_4^2+k_5^2+k_6^2)F_6.
\end{align*}
Using the above expressions it is easy to see that a number of terms in \eqref{4-harmonic-space-form} vanish and we obtain the following simplification
\begin{align*}
\nabla^7_TT+&K\nabla^5_TT
+K\langle\nabla^4_TT,T\rangle\nabla_TT
+K\langle\nabla^2_TT,T\rangle\nabla^3_TT=0,
\end{align*}
which, when expressed in terms of its Frenet frame, acquires the form
\begin{align*}
k_1&\bigg[
-(k_1^2+k_2^2)^3-k_2^2k_3^2(2k_1^2+2k_2^2+k_3^2+k_4^2)
+K\big((k_1^2+k_2^2)^2+k_2^2k_3^2\big)+2Kk_1^2(k_1^2+k_2^2)\bigg]F_2\\
&+k_1k_2k_3\bigg[(k_1^2+k_2^2)^2+(k_3^2+k_4^2)^2
+k_1^2k_3^2+2k_2^2k_3^2+k_4^2(k_1^2+k_2^2+k_5^2)
 \\
&-K(2k_1^2+k_2^2+k_3^2+k_4^2)
\bigg]F_4\\
&+k_1k_2k_3k_4k_5\bigg[
-(k_1^2+k_2^2+k_3^2+k_4^2+k_5^2+k_6^2)+K
\bigg]F_6=0.
\end{align*}
The claim now follows from testing this system with \(F_2,F_4,F_6\) completing the proof.
\end{proof}

\begin{proof}[Proof of Theorem \ref{thm:triharmonic-helices}]
The idea of the proof is to use the constraints \eqref{eq:triharmonic-curvatures}
and to perform a case by case analysis.
First of all we note that we have \(k_1\neq 0\) as we are considering a proper triharmonic curve.
\begin{enumerate}
\item If \(k_1\neq 0\) and \(k_i=0,i=1,2,3\) then
 we get \(k_1^2=2\) leading to the first class of curves.
 It was shown in \cite[Theorem 1.1]{MR4542687} that it actually 
 solves the equation for triharmonic curves.
 \item If \(k_1,k_2\neq 0\) and \(k_3=k_4=0\) we are in the situation detailed
 in Theorem \ref{thm:explicit-tri} leading to the second case.
 \item If \(k_1,k_2,k_3\neq 0\) and \(k_4=0\) the constraints 
 \eqref{eq:triharmonic-curvatures} acquire the form
 \begin{align}
 \label{constraint-k3}
  \sum_{i=1}^3k_i^2=1,\qquad (k_1^2+k_2^2)^2+k_2^2k_3^2&=2k_1^2+k_2^2.
 \end{align}
Using the first equation to eliminate \(k_3^2\) from the second one we find
\begin{align*}
k_1^2+k_2^2=2
\end{align*}
exploiting that \(k_1\neq 0\).
Reinserting this into the second equation of \eqref{constraint-k3} we find
\begin{align*}
k_2^2+k_2^2k_3^2=0
\end{align*}
leading to a contradiction such that this case cannot occur.
 
 \item If \(k_j\neq 0,j=1,\ldots 4\) the constraints \eqref{eq:triharmonic-curvatures} 
 are given by
 \begin{align}
 \label{constraint-k4}
  \sum_{i=1}^4k_i^2=1,\qquad (k_1^2+k_2^2)^2+k_2^2k_3^2&=2k_1^2+k_2^2.
 \end{align}
 Again, eliminating \(k_3^2\) from the second equation, making use of the first constraint,
 we find
 \begin{align*}
k_1^2(k_1^2+k_2^2)-k_2^2k_4^2=2k_1^2.
 \end{align*}
Using once more the first equation of \eqref{constraint-k4} to replace \(k_1^2+k_2^2\)
we arrive at
\begin{align*}
k_1^2(1+k_3^2+k_4^2)+k_2^2k_4^2=0
\end{align*}
leading to a contradiction again.
 
 \item If \(k_1\neq 0, k_2=0,k_3\neq 0,k_4=0\) 
the system \eqref{eq:triharmonic-curvatures} reduces to \(k_1^2=2\) leading to the
first claim of the theorem.
\item If \(k_1,k_2\neq 0,k_3=0,k_4\neq 0\) we get the condition \((k_1^2+k_2^2)^2=2k_1^2+k_2^2\) leading
to the second case of the theorem while
\(k_4\) can be arbitrary. However, it is a direct consequence of the Frenet equations \eqref{frenet-frame} that once \(k_3=0\) any geodesic curvature \(k_j,j\geq 4\) will no longer appear when 
expressing the equation for triharmonic curves in terms of its Frenet frame.
\end{enumerate}
The proof is now complete.
\end{proof}

A careful inspection of the proof of Theorem \ref{thm:tri-four-classification} shows that it is enough to know the structure
of the highest order derivatives of \(r\)-harmonic helices to 
obtain classification results if we 
assume that all geodesic curvatures of the curve are non-zero.
Hence, as a first step towards the proof of Theorem \ref{thm:structure}
we establish an expression for the iterated derivatives appearing
in the equation for \(r\)-harmonic curves \eqref{el-poly} suited to our particular analysis.

\begin{Lem}
Let \(\gamma\colon I\to M\) be an \(r\)-harmonic curve parametrized by arclength
whose geodesic curvature are all constant together with its Frenet frame \(\{F_j\},j=1,2r-2\).
\begin{enumerate}
 \item For \(2\leq l\leq 2r-3\) we have
 \begin{align}
 \label{derivative-general}
   \nabla^{2l-1}_TT=\sum_{j=1}^{l-2}a_jF_{2j}
   -\big(\prod_{i=1}^{2l-3}k_i\big)\big(\sum_{j=1}^{2l-2}k_j^2\big)F_{2l-2}
   +\big(\prod_{i=1}^{2l-1}k_i\big)F_{2l},
 \end{align}
 where \(a_j\) is a function of \(k_p,p=1,\ldots,l-2\).
 \item The highest derivative appearing in the equation for \(r\)-harmonic curves has the form
  \begin{align}
  \label{derivative-highest}
  \nabla^{2r-1}_TT=\sum_{j=1}^{r-2}b_jF_{2j}
   -\big(\prod_{i=1}^{2r-3}k_i\big)\big(\sum_{j=1}^{2r-2}k_j^2\big)F_{2r-2},
  \end{align}
  where \(b_j\) is a function of \(k_p,p=1,\ldots l-2\).
\end{enumerate}
\end{Lem}
\begin{proof}
The proof uses induction. 
Choosing \(l=3\) in \eqref{derivative-general} we get precisely the formula
derived in the proof of Theorem \ref{thm:tri-four-classification} confirming the base case.

For the induction step we differentiate \eqref{derivative-general} 
using the Frenet equations \eqref{frenet-frame} and find
\begin{align*}
\nabla^{2l}_TT=&\sum_{j=1}^{l-2}\tilde a_jF_{2j-1}+\sum_{j=1}^{l-2}a_jk_{2j}F_{2j+1} \\
&+\big(\prod_{i=1}^{2l-3}k_i\big)\big(\sum_{j=1}^{2l-2}k_j^2\big)k_{2l-3}F_{2l-3} 
-\big(\prod_{i=1}^{2l-3}k_i\big)\big(\sum_{j=1}^{2l-2}k_j^2\big)k_{2l-2}F_{2l-1}\\
&
-\big(\prod_{i=1}^{2l-1}k_i\big)k_{2l-1}F_{2l-1}
+\big(\prod_{i=1}^{2l-1}k_i\big)k_{2l}F_{2l+1},
\end{align*}
where \(\tilde a_j\) is again a function of the \(k_p,p=1,\ldots l-1\).
Differentiating again using \eqref{frenet-frame} then yields
\begin{align*}
\nabla^{2l+1}_TT=&-\sum_{j=1}^{l-2}\tilde a_jk_{2j-2}F_{2j-2}
+\sum_{j=1}^{l-2}\tilde a_jk_{2j-1}F_{2j}\\
&-\sum_{j=1}^{l-2}\tilde a_jk_{2j}^2F_{2j}
+\sum_{j=1}^{l-2}\tilde a_jk_{2j}k_{2j+1}F_{2j+2}\\
&-\big(\prod_{i=1}^{2l-3}k_i\big)\big(\sum_{j=1}^{2l-2}k_j^2\big)k_{2l-3}k_{2l-4}F_{2l-4}
+\big(\prod_{i=1}^{2l-3}k_i\big)\big(\sum_{j=1}^{2l-2}k_j^2\big)k^2_{2l-3}F_{2l-2} \\
&+\big(\prod_{i=1}^{2l-3}k_i\big)\big(\sum_{j=1}^{2l-2}k_j^2\big)k^2_{2l-2}F_{2l-2}
-\big(\prod_{i=1}^{2l-3}k_i\big)\big(\sum_{j=1}^{2l-2}k_j^2\big)k_{2l-2}k_{2l-1}F_{2l}\\
&+\big(\prod_{i=1}^{2l-1}k_i\big)k_{2l-1}k_{2l-2}F_{2l-2}
-\big(\prod_{i=1}^{2l-1}k_i\big)k_{2l-1}^2F_{2l} \\
&-\big(\prod_{i=1}^{2l-1}k_i\big)k_{2l}^2F_{2l}+
\big(\prod_{i=1}^{2l+1}k_i\big)F_{2l+2}.
\end{align*}
Now, it is straightforward to see that
\begin{align*}
\big(\prod_{i=1}^{2l-3}k_i\big)\big(\sum_{j=1}^{2l-2}k_j^2\big)k_{2l-2}k_{2l-1}F_{2l}
+\big(\prod_{j=1}^{2l-1}k_i\big)k_{2l-1}^2F_{2l}+\big(\prod_{i=1}^{2l-1}k_i\big)k_{2l}^2F_{2l}
=\big(\prod_{i=1}^{2l-1}k_i\big)\big(\sum_{j=1}^{2l}k_j^2\big)F_{2l}.
\end{align*}
Hence, we may conclude that 
\begin{align*}
\nabla^{2l+1}_TT=\sum_{j=1}^{l-1}\tilde a_jF_{2j}
-\big(\prod_{i=1}^{2l-1}k_i\big)\big(\sum_{j=1}^{2l}k_j^2\big)F_{2l}
+\big(\prod_{i=1}^{2l+1}k_i\big)F_{2l+2}
 \end{align*}
completing the induction step and thus establishing the first claim.
The second formula follows from the first one taking into account
that for an \(r\)-harmonic curve we have \(k_{2r-1}=0\).

\end{proof}

We are now ready to give the proof of Theorem \ref{thm:structure}.

\begin{proof}[Proof of Theorem \ref{thm:structure}]
First, we rewrite the equation for \(r\)-harmonic curves
in space forms
\eqref{polyharmonic-spaceform},
extracting the two leading derivatives,
as follows
\begin{align}
\label{euler-lagrange-split}
\nabla^{2r-1}_TT&+K\nabla^{2r-3}_TT
-K\langle T,\nabla^{2r-3}_TT\rangle T \\
\nonumber&+\sum_{l=1}^{r-2}(-1)^l\big(\langle T,\nabla^l_TT\rangle\nabla_T^{2r-3-l}T
-\langle T,\nabla_T^{2r-3-l}T\rangle\nabla_T^lT\big)=0.
\end{align}
Testing this equation with \(F_{2r-2}\) we obtain 
\begin{align*}
\langle& \nabla^{2r-1}_TT,F_{2r-2}\rangle
+K\langle\nabla^{2r-3}_TT,F_{2r-2}\rangle \\
&+\sum_{l=1}^{r-2}(-1)^l\big(\langle T,\nabla^l_TT\rangle
\langle\nabla_T^{2r-3-l}T,F_{2r-2}\rangle
-\langle T,\nabla_T^{2r-3-l}T\rangle
\langle\nabla_T^lT,F_{2r-2}\rangle\big)=0.
\end{align*}
Regarding the last two terms, we make the following splitting
into even and odd addends
\begin{align*}
\sum_{l=1}^{r-2}(-1)^l\langle\nabla_T^{2r-3-l}T,F_{2r-2}\rangle= 
\sum_{l=1}^{\frac{r-2}{2}}\big(\langle\nabla_T^{2r-3-2l}T,F_{2r-2}\rangle
-\langle\nabla_T^{2r-3-(2l-1)}T,F_{2r-2}\rangle\big).
\end{align*}
As \(2r-3-(2l-1)=2r-2-2l\) is clearly even it is a direct consequence of the 
Frenet equations \eqref{frenet-frame} that \(\nabla_T^{2r-3-(2l-1)}T\)
can be written as
\begin{align*}
\nabla_T^{2r-3-(2l-1)}T=\sum_ja_jF_{2j+1},
\end{align*}
where \(a_j\) are functions of the geodesic curvatures \(k_j\),
such that 
\begin{align*}
\langle\nabla_T^{2r-3-(2l+1)}T,F_{2r-2}\rangle=0.
\end{align*}
Secondly, using \eqref{derivative-general}, we find
\begin{align*}
 \langle\nabla_T^{2r-3-2l}T,F_{2r-2}\rangle=&
   \sum_{j=1}^{r-l-3}a_j\underbrace{\langle F_{2j},F_{2r-2}\rangle}_{=0} \\
   &-\big(\prod_{i=1}^{2r-2l-5}k_i\big)\big(\sum_{j=1}^{2r-2l-4}k_j^2\big)
   \underbrace{\langle F_{2(r-l-2)},F_{2r-2}\rangle}_{=0} \\
   &+\big(\prod_{i=1}^{2r-2l-3}k_i\big)\underbrace{\langle F_{2(r-l-1)},F_{2r-2}}_{=0}\rangle
 \end{align*}
for all \(l\geq 1\).
Hence, we may conclude that
\begin{align*}
\sum_{l=1}^{r-2}(-1)^l\langle\nabla_T^{2r-3-l}T,F_{2r-2}\rangle=0
\end{align*}
and by the same reasoning we can also deduce that
\begin{align*}
\langle\nabla_T^lT,F_{2r-2}\rangle=0,\qquad 1\leq l\leq r-2.
\end{align*}
Now, using \eqref{derivative-general} and \eqref{derivative-highest}
we obtain from the equation for \(r\)-harmonic curves \eqref{euler-lagrange-split} that
\begin{align*}
0=&\langle\nabla^{2r-1}_TT,F_{2r-2}\rangle
+K\langle\nabla^{2r-3}_TT,F_{2r-2}\rangle\\
=&-\big(\prod_{i=1}^{2r-3}k_i\big)\big(\sum_{j=1}^{2r-2}k_j^2\big)
+K\big(\prod_{i=1}^{2r-3}k_i\big).
\end{align*}
This completes the proof. 
\end{proof}

\subsection{Evidence in support of Conjecture 1}
\label{subsection-conjecture}
In this subsection we will collect a number of results which support 
the statement of Conjecture \ref{conjecture}.

First, recall that the equation for a triharmonic curve on a general Riemannian
manifold is given by
\begin{align}
\label{triharmonic-general}
0=\tau_3(\gamma)=\nabla^5_TT+R^M(\nabla^3_TT,T)T-R^M(\nabla^2_TT,\nabla_TT)T,
\end{align}
which is precisely \eqref{el-poly} for \(r=3\).

\begin{Prop}
\label{prop-conservation-triharmonic}
Let \(\gamma\colon I\to M\) be a triharmonic curve
parametrized with respect to arclength.
Then the following equation holds
\begin{align}
\label{conservation-triharmonic}
\frac{d^3}{ds^3}|\nabla_TT|^2-\frac{d}{ds}|\nabla^2_TT|^2=&0.
\end{align}
\end{Prop}

\begin{proof}
To obtain the first identity 
we multiply \eqref{triharmonic-general} by \(T\) and calculate
\begin{align*}
0=&\langle\nabla^5_TT,T\rangle \\
=&\frac{d}{ds}\langle\nabla^4_TT,T\rangle-\langle\nabla^4_TT,\nabla_TT\rangle \\
=&\frac{d^2}{ds^2}\langle\nabla^3_TT,T\rangle-2\frac{d}{ds}\langle\nabla^3_TT,\nabla_TT\rangle
+\langle\nabla^3_TT,\nabla^2_TT\rangle \\
=&\frac{d^3}{ds^3}\langle\nabla^2_TT,T\rangle-3\frac{d^2}{ds^2}\langle\nabla^2_TT,\nabla_TT\rangle
+\frac{5}{2}\frac{d}{ds}|\nabla^2_TT|^2.
\end{align*}
To finish the proof we note that \(\langle\nabla^2_TT,T\rangle=-|\nabla_TT|^2\)
which holds as \(\gamma\) is parametrized with respect to arclength.
\end{proof}

\begin{Cor}
Let \(\gamma\colon I\to N\) be a triharmonic curve
parametrized with respect to arclength.
Then the following conservation law holds
\begin{align}
\label{conservation-triharmonic-constant}
\frac{d^2}{ds^2}|\nabla_TT|^2-|\nabla^2_TT|^2=&c_1
\end{align}
for some \(c_1\in\R\).
\end{Cor}
\begin{proof}
This is a direct consequence of the conservation law \eqref{conservation-triharmonic}.
\end{proof}

Choosing a Frenet frame along \(\gamma\) equation \eqref{conservation-triharmonic-constant} implies
\begin{align*}
(k_1')^2+2k_1k_1''-k_1^4-k_1^2k_2^2=c_1,
\end{align*}
where \(k_i,i=1,2\) represent the curvatures of the curve \(\gamma\).
In the case of \(c_1=0\) this equation is solved by 
\begin{align*}
k_1=\frac{\alpha}{s},\qquad k_2=\frac{\beta}{s},\qquad \alpha^2+\beta^2=5,
\end{align*}
which gives rise to the triharmonic curve with non-constant geodesic 
curvature and torsion constructed in 
\cite{MR4308322}.

In the following, we extend the previous analysis to the case of \(4\)-harmonic curves.
These are solutions of 
\begin{align}
\label{euler-lagrange-4}
0=\tau_4(\gamma)=\nabla^7_TT+R^M(\nabla^5_TT,T)T
-R^M(\nabla^4_TT,\nabla_TT)T
+R^M(\nabla^3_TT,\nabla^2_TT)T,
\end{align}
which is precisely \eqref{el-poly} in the case of \(r=4\).

\begin{Prop}
Let \(\gamma\colon I\to M\) be a 4-harmonic curve
parametrized with respect to arclength.
Then the following conservation law holds
\begin{align}
\label{conservation-4harmonic}
\frac{d^5}{ds^5}|\nabla_TT|^2
-2\frac{d^3}{ds^3}|\nabla^2_TT|^2
+\frac{d}{ds}|\nabla^3_TT|^2=0.
\end{align}

\end{Prop}
\begin{proof}
Testing \eqref{euler-lagrange-4} with \(T\) 
a direct calculation yields the following identity
\begin{align*}
0=\frac{d^2}{ds^2}\langle\nabla^5_TT,T\rangle
-2\frac{d}{ds}\langle\nabla^5_TT,\nabla_TT\rangle
+\langle\nabla^5_TT,\nabla^2_TT\rangle.
\end{align*}

From the proof of Proposition \ref{prop-conservation-triharmonic} we know that
\begin{align*}
\langle\nabla^5_TT,T\rangle=-\frac{5}{2}\frac{d^3}{ds^3}|\nabla_TT|^2
+\frac{5}{2}\frac{d}{ds}|\nabla^2_TT|^2.
\end{align*}

Moreover, we have
\begin{align*}
\langle\nabla^5_TT,\nabla_TT\rangle=
\frac{d^4}{ds^4}\frac{1}{2}|\nabla_TT|^2
-\frac{d^2}{ds^2}2|\nabla^2_TT|^2+|\nabla^3_TT|^2.
\end{align*}

Finally, a direct calculation shows
\begin{align*}
\langle\nabla^5_TT,\nabla^2_TT\rangle=
\frac{1}{2}\frac{d^3}{ds^3}|\nabla^2_TT|^2
-\frac{3}{2}\frac{d}{ds}|\nabla^3_TT|^2.
\end{align*}
The claim follows by combing the different equations.
\end{proof}

A dimensional analysis of the conservation law \eqref{conservation-4harmonic} suggests the following:
Assume that we are looking for a 4-harmonic curve with non-constant
geodesic curvature \(k_1\). Inspecting the terms in \eqref{conservation-4harmonic} suggests that \(k_1=\frac{C}{s^2}\)
as all three terms scale as \(\frac{1}{s^9}\).
Hence, one may expect that there exist 4-harmonic curves with non-constant geodesic curvature 
\(k_1=\frac{\alpha}{s^2},\alpha\in\R\).

Finally, we note that the construction of conservation laws can
be carried out for polyharmonic curves as well by multiplying \eqref{el-poly} with \(T\)
and manipulating the resulting equation as we have demonstrated
for triharmonic and \(4\)-harmonic curves.
Again, a simple dimensional analysis then leads to the conclusion of Conjecture \ref{conjecture}.

\bibliographystyle{plain}
\bibliography{mybib}
\end{document}